\documentclass[12pt]{amsart}
\usepackage{amsmath,amscd,amssymb,amsfonts}
\newcommand{\ver}{Sept. 17, 2005, v.6}
\newcommand{\ssbull}{\raise.2ex\hbox{${\scriptscriptstyle\bullet}$}}
\newcommand{\scirc}{\raise.2ex\hbox{${\scriptstyle\circ}$}}
\newcommand{\mtim}{\hbox{$\times$}}
\newcommand{\msum}{\hbox{$\sum$}}
\newcommand{\mprod}{\hbox{$\prod$}}
\newcommand{\mopls}{\hbox{$\bigoplus$}}
\newcommand{\mcap}{\hbox{$\bigcap$}}
\newcommand{\mcup}{\hbox{$\bigcup$}}
\newcommand{\vep}{{\varepsilon}}
\newcommand{\bA}{{\mathbb A}}
\newcommand{\bC}{{\mathbb C}}
\newcommand{\bN}{{\mathbb N}}
\newcommand{\bQ}{{\mathbb Q}}
\newcommand{\bR}{{\mathbb R}}
\newcommand{\bZ}{{\mathbb Z}}
\newcommand{\boR}{{\mathbf R}}
\newcommand{\bone}{\hbox{\bf 1}}
\newcommand{\ba}{{\mathbf a}}
\newcommand{\bs}{{\mathbf s}}

\newcommand{\cD}{{\mathcal D}}
\newcommand{\cH}{{\mathcal H}}
\newcommand{\cJ}{{\mathcal J}}
\newcommand{\cM}{{\mathcal M}}
\newcommand{\cO}{{\mathcal O}}

\newcommand{\fa}{{\mathfrak a}}
\newcommand{\oM}{{\overline{M}}}
\newcommand{\tpr}{{\widetilde{pr}}}
\newcommand{\tM}{{\widetilde{M}}}
\newcommand{\tm}{{\widetilde{m}}}
\newcommand{\tilt}{{\widetilde{t}}}
\newcommand{\ttau}{{\widetilde{\tau}}}
\newcommand{\tx}{{\widetilde{x}}}
\newcommand{\tX}{{\widetilde{X}}}
\newcommand{\tZ}{{\widetilde{Z}}}
\newcommand{\tom}{{\widetilde{\omega}}}
\newcommand{\SP}{\hbox{\rm Sp}}
\newcommand{\Spec}{\hbox{\rm Spec}}
\newcommand{\codim}{\hbox{\rm codim}}
\newcommand{\Gr}{\hbox{\rm Gr}}
\newcommand{\Ker}{\hbox{\rm Ker}}
\renewcommand{\Re}{\hbox{\rm Re}}
\renewcommand{\Im}{\hbox{\rm Im}}
\newcommand{\DR}{\hbox{\rm DR}}
\newcommand{\pr}{\hbox{\it pr}}
\newcommand{\reg}{\text{\rm reg}}
\newcommand{\Sing}{\hbox{\rm Sing}}

\newcommand{\simto}{\buildrel\sim\over\longrightarrow}

\begin{document}
\title[Bernstein-Sato Polynomials of Arbitrary Varieties]
{Bernstein-Sato Polynomials of Arbitrary Varieties}
\author{Nero Budur}
\address{Department of Mathematics, The Johns Hopkins
University, Baltimore, MD 21218 USA}
\email{nbudur@math.jhu.edu}
\author{Mircea Musta\c{t}\v{a}}
\address{Department of Mathematics, University of Michigan,
Ann Arbor, MI 48109, USA}
\email{mmustata@umich.edu}
\author{Morihiko Saito}
\address{RIMS Kyoto University, Kyoto 606-8502 Japan}
\email{msaito@kurims.kyoto-u.ac.jp}
\date{\ver}
\begin{abstract}
We introduce the notion of Bernstein-Sato polynomial of
an arbitrary variety (which is not necessarily reduced
nor irreducible), using the theory of
$ V $-filtrations of M.~Kashiwara and B.~Malgrange.
We prove that the decreasing filtration by multiplier ideals
coincides essentially with the restriction of the
$ V $-filtration.
This implies a relation between the roots of the Bernstein-Sato
polynomial and the jumping coefficients of the multiplier
ideals, and also a criterion for rational singularities in
terms of the maximal root of the polynomial
in the case of a reduced complete intersection.
These are generalizations of the hypersurface case.
We can
calculate the polynomials explicitly in the case of monomial
ideals.
\end{abstract}
\maketitle

\bigskip
\centerline{\bf Introduction}

\bigskip\noindent
The notion of Bernstein-Sato polynomial (i.e.
$ b $-function) for a function was introduced
independently by J.N.~Bernstein [3] and M.~Sato [33].
This theory was then developed by J.-E.~Bj\"ork [4],
B.~Malgrange [21], and by members of the Sato school:
M.~Kashiwara [15] and T.~Yano [37] among others.
Related to the theory of vanishing cycles of P.~Deligne [10],
it further led to the theory of
$ V $-filtrations of Kashiwara [16] and Malgrange [22].
It is well-known and is easy to show that the
$ b $-function depends only on the hypersurface defined by
the function.
Motivated by our previous work [6] related to multiplier ideals [20],
we had the necessity to generalize the notion of
$ b $-function to the case of arbitrary subvarieties.

Let
$ Z $ be a (not necessarily reduced nor irreducible) complex
algebraic variety
embedded in a smooth affine variety
$ X $.
Let
$ f_{1},\dots, f_{r} $ be nonzero generators of the ideal of
$ Z $ (i.e.
$ f_{j} \ne 0 $).
Let
$ \cD_{X} $ be the sheaf of linear differential operators on $ X $.
It acts naturally on
$ \cO_{X}[\mprod_{i}f_{i}^{-1},s_{1},\dots, s_{r}]
\mprod_{i}f_{i}^{{s}_{i}} $,
where the
$ s_{i} $ are independent variables.
We define a
$ \cD_{X} $-linear action of
$ t_{j} $ on it by
$ t_{j}(s_{i}) = s_{i} + 1 $ if
$ i = j $, and
$ t_{j}(s_{i}) = s_{i} $ otherwise.
In particular,
$ t_{j}\mprod_{i}f_{i}^{{s}_{i}} = f_{j}\mprod_{i}f_{i}^{{s}_{i}} $,
and the action of
$ t_{j} $ is bijective.
Let
$ s_{i,j} = s_{i}t_{i}^{-1}t_{j} $,
and
$ s = \msum_{i}s_{i} $.
The Bernstein-Sato polynomial (i.e. the
$ b $-function)
$ b_{f}(s) $ of
$ f := (f_{1}, \dots, f_{r}) $ is defined to be the monic
polynomial of the lowest degree in
$ s $ satisfying the relation
$$
b_{f}(s)\mprod_{i}f_{i}^{{s}_{i}} =
\msum_{k=1}^r P_{k}t_{k}
\mprod_{i}f_{i}^{{s}_{i}},
$$
where the
$ P_{k} $ belong to the ring generated by
$ \cD_{X} $ and the
$ s_{i,j} $;
see (2.4) and (2.10) for other formulations.
Note that if we require instead $ P_k\in\cD_{X}[s_1,\ldots,s_r] $,
then there are examples in which there is no nonzero such $b_{f}(s)$
(see (4.5)).
The above definition gives a natural generalization of the
$ b $-function in the hypersurface case (i.e.
$ r = 1 $, see [3]), and this normalization of the
$ b $-function is the same as in [6], [15], [31].
Since our definition of the
$ b $-function is closely related to the
$ V $-filtration of Kashiwara [16] and Malgrange [22],
its existence and the rationality of its roots
follow easily from their theory (see also [13]).
We can show moreover that the denominators of the roots of the
$ b $-function are related to the multiplicities of the
divisor obtained by an embedded resolution of
$ (X,Z) $, see (3.10).

By (2.4) below this
$ b_{f}(s) $ coincides (up to a shift of variable) with the polynomial
$ b_{L} $ of the minimal degree satisfying the relation in
[27], I, 3.1.1 in the algebraic setting (see also [13], 2.13), if
$ L(s) = s_{1}+\cdots + s_{r} $ and if
$ \cM, U $ in loc.~cit. are chosen appropriately, see (5.2) below.
This polynomial was used to prove the main theorem in [27], I.
However,
$ b_{f}(s) $ is slightly different from a polynomial in [27], II,
Prop.~1.1,
because our definition requires certain additional binomial
polynomials as in (2.10) below, and Theorem~2 does not hold for that
polynomial without the additional term, see (5.4) below.
In [27] Sabbah proved the existence of nonzero polynomials of
{\it several variables}
which satisfy a functional equation similar to the above one,
see also [1], [8], [13].
However, its relation with
$ b_{f}(s) $ seems to be quite nontrivial, see (5.1) below.

Let
$ b_{Z}(s) = b_{f}(s-r') $ with
$ r' = \codim_{X}Z $.
This normalization is the same as in [3].
Using the
$ V $-filtration of Kashiwara and Malgrange, we can show that our
$ b_{f}(s) $ is independent of the choice of a system of generators
$ f = (f_{1},\dots,f_{r}) $ if
$ Z $ and
$ \dim X $ are fixed (but
$ r $ can vary), and that
$ b_{Z}(s) $ depends only on
$ Z $, see (2.5).
This was rather surprising because the assertion does not hold
for the polynomial in [27], II, 1.1.
If
$ Z $ is not affine, the
$ b $-function can be defined to be the least common
multiple of the local
$ b $-functions (shifted appropriately if
$ Z $ is not equidimensional).
For
$ g \in \cO_{X} $,
we similarly define the
$ b $-function
$ b_{f,g}(s) $ with
$ \mprod_{i}f_{i}^{{s}_{i}} $ replaced by
$ g\mprod_{i}f_{i}^{{s}_{i}} $.

Our first main theorem is concerning the multiplier ideals and
the $V$-filtration.
For a positive rational number
$ \alpha $, the multiplier ideal sheaf
$ \cJ(X,\alpha Z) $ is defined by taking an embedded
resolution of
$ (X,Z) $ (see [20] and also (3.1) below)
or using the local integrability of
$ |g|^{2}/(\msum_{i}|f_{i}^{2}|)^{\alpha} $ for
$ g \in \cO_{X} $ (see [25]).
This gives a decreasing sequence of ideals, and
there are positive rational numbers
$ 0 < \alpha_{1} < \alpha_{2} < \cdots $ such that
$ \cJ(X,\alpha_{j}Z) = \cJ(X,\alpha Z) \ne \cJ(X,\alpha_{j+1}Z) $
for
$ \alpha_{j} \le \alpha < \alpha_{j+1} \,(j \ge 0) $ where
$ \alpha_{0} = 0 $ and
$ \cJ(X,\alpha_{0}Z) = \cO_{X} $.
These
$ \alpha_{j} $ for
$ j > 0 $ are called the jumping coefficients,
and the minimal jumping coefficient
$ \alpha_{1} $ is called the log-canonical threshold.

We will denote by
$ V $ the filtration on
$ \cO_{X} $ induced by the
$ V $-filtration of Kashiwara [16] and Malgrange [22] along
$ Z $.

\medskip\noindent
{\bf Theorem~1.}
{\it We have
$ V^{\alpha}\cO_{X} = \cJ(X,\alpha Z) $ if
$ \alpha $ is not a jumping coefficient.
In general,
$ \cJ(X,\alpha Z) = V^{\alpha +\vep}\cO_{X} $ and
$ V^{\alpha}\cO_{X} = \cJ(X,(\alpha -\vep)Z) $ for any
$ \alpha \in \bQ $, where
$ \vep > 0 $ is sufficiently small.
}

\medskip
This implies the following relation between the roots of the
$ b $-function and the jumping coefficients.

\medskip\noindent
{\bf Theorem~2.}
{\it The log-canonical threshold of
$ (X,Z) $ coincides with the smallest root
$ \alpha'_{f} $ of
$ b_{f}(-s) $ {\rm (}in particular,
$ \alpha'_{f} > 0 )$, and any jumping coefficients of
$ (X,Z) $ in
$ [\alpha'_{f},\alpha'_{f}+1) $ are roots of
$ b_{f}(-s) $.
}

\medskip
In the case $Z$ is a reduced complete intersection,
we have an analogue of the adjoint ideal of a reduced divisor
(see [11], [23], [35] and also [6], 3.7),
and the maximal root
$ -\alpha'_{f} $ of
$ b_{f}(s) $ can be used for a criterion of rational
singularities as follows.

\medskip\noindent
{\bf Theorem~3.}
{\it Assume
$ Z $ is a reduced complete intersection of
codimension
$ r $ in
$ X $.
Let
$ \pi : \tZ \to Z $ be a resolution of singularities, and
set
$ \tom_{Z} = \pi_{*}\omega_{\tZ} $ where
$ \omega_{\tZ} $ denotes the dualizing sheaf.
Then there is a coherent ideal
$ \cJ(X,rZ)' $ of
$ \cO_{X} $ such that
$ \cJ(X,rZ) \subset \cJ(X,rZ)' \subset \cJ(X,(r-\vep)Z) $
for
$ 0 < \vep \ll 1 $ and
$$
\omega_{X}\otimes (\cO_{X}/\cJ(X,rZ)') = \omega_{Z}/\tom_{Z}.
$$
In particular,
$ \omega_{X}\otimes (\cJ(X,(\alpha-\vep)Z)/
\cJ(X,\alpha Z)) $ is isomorphic to a subquotient of
$ \omega_{Z}/\tom_{Z} $ for
$ \alpha < r $ and
$ 0 < \vep \ll 1 $.
}

\medskip\noindent
{\bf Theorem~4.}
{\it With the assumption of Theorem~{\rm 3},
$ Z $ has at most rational singularities if and only if
$ \alpha'_{f} = r $ and its multiplicity is
$ 1 $.
}

\medskip
The above theorems generalize the corresponding
results in the hypersurface case, see [6],
[12], [15], [18], [31], [34].
In Theorem~4,
$ r $ is always a root of
$ b_{f}(-s) $ by restricting
to the smooth points of
$ Z $.
Let
$ \alpha_{f} $ denote the minimal root of
$ b_{f}(-s)/(-s+r) $.
Then the criterion for rational singularities in Theorem~4
is equivalent to the condition
$ \alpha_{f} > r $.
In the hypersurface isolated singularity case, it is known that
$ \alpha_{f} $ coincides with the Arnold exponent.

In the case of rational singularities, Theorem~3 implies

\medskip\noindent
{\bf Corollary~1.}
{\it Under the assumption of Theorem~{\rm 3}, assume
further that
$ Z $ has at most rational singularities, or more generally,
$ \alpha'_{f} = r $.
Then the jumping coefficients are the integers
$ \ge r $, and
$ \cJ(X,jZ) = I_{Z}^{j-r+1} $ for
$ j \ge r $ where
$ I_{Z} $ is the ideal sheaf of
$ Z $.
}

\medskip
From Theorem~1 (together with (2.1.3) below) we deduce also
the following description of multiplier ideals in
terms of $b$-functions.

\medskip\noindent
{\bf Corollary~2.}
{\it For
$ g \in \cO_{X} $, it belongs to
$ \cJ(X,\alpha Z) $ if and only if all the roots of
$ b_{f,g}(-s) $ are strictly greater than
$ \alpha $.
}

\medskip
In the hypersurface case this is due to [26] if
$ \cJ(X,\alpha Z) $ is replaced by
$ V^{>\alpha}\cO_{X} $.
If the ideal of
$ Z $ is a monomial ideal, we can calculate the
$ b $-function in some cases (see Section 4).
Other examples can be deduced from the following
Thom-Sebastiani type theorem,
which is compatible with the similar theorem for
multiplier ideals (see [24]) via Theorem~2. Note that
it is different from a usual Thom-Sebastiani type theorem,
which applies to the sum of
the pull-backs of two functions, see e.g. [37].

\medskip\noindent
{\bf Theorem~5.}
{\it For
$ f : X \to \bA^{r} $ and
$ g : Y \to \bA^{r'} $,
let
$ h = f\times g : X\times Y \to \bA^{r+r'} $.
Write
$ b_{f}(s) = \prod_{\alpha} (s+\alpha )^{n_{\alpha}},
b_{g}(s) = \prod _{\beta} (s+\beta )^{m_{\beta}},
b_{h}(s) = \prod _{\gamma} (s+\gamma )^{q_{\gamma}} $.
Then
$ q_{\gamma} = \max\{n_{\alpha} + m_{\beta} - 1 \,|\,
n_{\alpha} > 0, m_{\beta} > 0, \alpha + \beta = \gamma \} $.
}

\medskip
This paper is organized as follows.
In Sect.~1, we review the theories of
$ V $-filtrations and specializations due to Kashiwara, Malgrange
and Verdier.
In Sect.~2, we define the Bernstein-Sato polynomial for an
arbitrary variety (which is not necessarily reduced
nor irreducible), and prove its existence and
well-definedness together with Theorem~5.
In Sect.~3, we show the bistrictness of the direct image,
and prove Theorems~1--4.
In Sect.~4, we treat the monomial ideal case, and calculate
some examples.
In Sect.~5, we explain the relation with the Bernstein-Sato
polynomials in other papers.

\medskip\noindent
{\bf Convention.}
In this paper, a variety means a (not necessarily reduced nor
irreducible) separated scheme of finite type over
$ \bC $, and a point of a variety means a closed point.
In particular, the underlying set of a variety coincides
with that of the associated analytic space.

\bigskip\bigskip
\centerline{{\bf 1. Filtration of Kashiwara and Malgrange}}

\bigskip\noindent
In this section we review the theories of
$ V $-filtrations and specializations due to Kashiwara, Malgrange
and Verdier.

\medskip\noindent
{\bf 1.1.
$ V $-Filtration.}
Let
$ Z $ be a smooth closed subvariety of codimension
$ r $ in a smooth complex algebraic variety
$ X $, and
$ I_{Z} $ be the ideal sheaf of
$ Z $ in
$ X $.
Let
$ \cD_{X} $ be the sheaf of linear differential operators on
$ X $.
The filtration
$ V $ on
$ \cD_{X} $ along
$ Z $ is defined by
$$
V^{i}\cD_{X} = \{P \in \cD_{X} : PI_{Z}^{j} \subset
I_{Z}^{j+i}\quad \text{for any}\,\, j \ge 0 \},
$$
where
$ I_{Z}^{i} = \cO_{X} $ for
$ i \le 0 $.
Let
$ (x_{1}, \dots, x_{n}) $ be a local coordinate system (i.e.
inducing an etale morphism to
$ \bA^{n}) $ such that
$ Z = \{x_{i} = 0 \,(i\le r)\} $.
Let
$ \partial_{x_{i}} = \partial /\partial x_{i} $.
Then
$ V^{k}\cD_{X} $ is generated over
$ \cO_{X} $ by
$$
\mprod_{j\le r}x_{j}^{\mu_{j}}\mprod_{i\le n}
\partial_{x_{i}}^{{\nu}_{i}}\quad \text{with}\quad
\msum_{j\le r}\mu_{j} - \msum_{i\le r}\nu_{i} = k,
\leqno(1.1.1)
$$
where
$ (\mu_{1}, \dots, \mu_{r}) \in \bN^{r} $ and
$ (\nu_{1}, \dots, \nu_{n}) \in \bN^{n} $.

We say that a filtration
$ V $ is discretely indexed by
$ \bQ $ if there is an increasing sequence of rational numbers
$ \{\alpha_{j}\}_{j\in \bZ} $ such that
$ \lim_{j\to-\infty}\alpha_{j} = -\infty $,
$ \lim_{j\to+\infty}\alpha_{j} = +\infty $, and
$ V^{\alpha} $ for
$ \alpha \in (\alpha_{j},\alpha_{j+1}) $ depends only on
$ j $.
We say that a decreasing filtration
$ V $ is left-continuously indexed if
$ V^{\alpha}M = \mcap_{\beta <\alpha} V^{\beta}M $ for any
$ \alpha $.
Let
$ \theta $ be a (locally defined) vector field such that
$ \theta \in V^{0}\cD_{X} $ and whose action on
$ I_{Z}/I_{Z}^{2} $ is the identity.
For a
$ \cD_{X} $-module
$ M $,
the
$ V $-filtration of Kashiwara [16] and Malgrange [22] along
$ Z $ is an exhaustive decreasing filtration
$ V $ which is indexed discretely and left-continuously by
$ \bQ $,
and satisfies the following conditions (see also [19], [26]):

\medskip\noindent
(i) The
$ V^{\alpha}M $ are coherent
$ V^{0}\cD_{X} $-submodules of
$ M $.

\noindent
(ii)
$ V^{i}\cD_{X}V^{\alpha}M \subset V^{\alpha +i}M $
for any
$ i \in \bZ, \alpha \in \bQ $.

\noindent
(iii)
$ V^{i}\cD_{X}V^{\alpha}M = V^{\alpha +i}M $ for any
$ i > 0 $, if
$ \alpha \gg 0 $.

\noindent
(iv)
$ \theta - \alpha + r $ is nilpotent on
$ \Gr_{V}^{\alpha}M $.

\medskip\noindent
Here
$ \Gr_{V}^{\alpha}M := V^{\alpha}M/V^{>\alpha}M $ with
$ V^{>\alpha}M = \mcup_{\beta >\alpha} V^{\beta}M $.
Condition (iii) is equivalent to the condition:
$ V^{1}\cD_{X}V^{\alpha}M = V^{\alpha +1}M $ for
$ \alpha \gg 0 $, assuming condition (ii).
Condition (iv) is independent of the choice of
$ \theta $.
The shift of index by
$ r $ in condition (iv) is necessary to show an assertion
related to the independence of embeddings in smooth varieties.

By the theory of Kashiwara [16] and
Malgrange [22], there exists uniquely the
$ V $-filtration on $M$ indexed by
$ \bQ $ if
$ M $ is regular holonomic and quasi-unipotent, see also
Remark (1.2)(iv) below.
(A holonomic
$ \cD_{X} $-module is called regular and quasi-unipotent, if
its pull-back by any morphism of a curve to
$ X $ has regular singularities and quasi-unipotent local
monodromies.)

For example, if
$Z$ is smooth, then the $V$-filtration on $M=\cO_{X}$ along $Z$
is given by $V^{\alpha}\cO_{X}=I_Z^{\lceil\alpha\rceil-r}$, where
$\lceil\alpha\rceil$ is the smallest integer $\geq\alpha$.

\medskip\noindent
{\bf 1.2.~Remarks.} (i)
The above conditions are enough to characterize the
$ V $-filtration uniquely.
Indeed, if there are two filtrations
$ V, V' $ satisfying the above conditions, we have
by condition (iv)
$$
\Gr_{V}^{\alpha}\Gr_{V'}^{\beta}M = 0\quad\text{for}\,\,
\alpha \ne \beta.
\leqno(1.2.1)
$$
So it is enough to show that
$$
V^{\beta +k} \subset V'{}^{\beta} \subset V^{\beta -k}
\quad\text{for}\,\, k \gg 0,
\leqno(1.2.2)
$$
because it implies that
$ V $ and
$ V' $ induce finite filtrations on
$ V^{\alpha}/V^{\beta} $,
$ V'{}^{\alpha}/V'{}^{\beta} $ for any
$ \alpha < \beta $ and these induced filtrations coincide
by (1.2.1).
The second inclusion of (1.2.2) follows from condition (i),
and we need condition (iii) to show the first inclusion.

(ii) Conditions (i)--(iv) in (1.1) imply that equality holds
in condition (ii) if
$ i <0 $ and
$ \alpha \ll 0 $.
Indeed,
$ V^{\beta}M $ generates
$ M $ over
$ \cD_{X} $ for
$ \beta $ sufficiently small (considering an increasing sequence
of
$ \cD_{X} $-submodules generated by
$ V^{\beta}M) $, and
$ V $ coincides with the filtration defined by
$ V'{}^{\alpha}M = V^{\alpha}M $ for
$ \alpha \ge \beta $, and
$ V'{}^{\alpha}M = V^{-i}\cD_{X}V^{\alpha+i}M $ otherwise,
where
$ i $ is an integer such that
$ \beta \le \alpha+i < \beta + 1 $.

(iii) Let
$ M' $ be a
$ \cD_{X} $-submodule of
$ M $.
Then the restriction of
$ V $ to
$ M' $ satisfies the conditions of the
$ V $-filtration.
Indeed,
$ \mopls_{i\in\bN}V^{\alpha+i}M' $ is finitely generated over
$ \mopls_{i\in\bN}V^{i}\cD_{X} $ by the noetherian property, because
$ \mopls_{i,p\in\bN}\Gr_{p}^{F}V^{i}\cD_{X} =
\mopls_{i,p\in\bN}I_{Z}^{i+p}\Gr_{p}^{F}\cD_{X} $
by (1.1.1), where
$ F $ is the filtration by the order of differential operators,
see [26].

(iv) If we do not assume that
$ M $ is quasi-unipotent, then, after choosing an order on $\bC$
(for example, such that
$ \alpha<\beta $ if and only if
$ \Re\,\alpha < \Re\,\beta $ or
$ \Re\,\alpha = \Re\,\beta $ and
$ \Im\,\alpha < \Im\,\beta $),
there is a $V$-filtration indexed by $\bC$, see also [27].
If
$ M $ is quasi-unipotent, we see that
$ V $ is actually indexed by
$ \bQ $ using (1.3.1) below, because the assertion
in the codimension one case is well-known
(and is easily proved by using a resolution of singularities).

\medskip\noindent
{\bf 1.3.~ Specialization.}
With the notation of (1.1), let
$$
\tX = {\mathcal Spec}_{X}(\mopls_{i\in \bZ} I_{Z}^{-i}\otimes
t^{i}),
$$
where
$ I_{Z}^{-i} = \cO_{X} $ for
$ i \ge 0 $.
This is an open subvariety of the blow-up of
$ X\mtim \bC $ along
$ Z\mtim \{0\} $.
There is a natural morphism
$ p : \tX \to \bA^{1} := \Spec\, \bC[t] $ whose fiber
over
$ 0 $ is the tangent cone
$$
T_{Z}X = {\mathcal Spec}_{X}(\mopls_{i\le 0} I_{Z}^{-i}/I_{Z}^{-i+1}
\otimes t^{i}),
$$
and
$ \tX^{*} := p^{-1}(\bA^{1}\backslash \{0\}) $ is
isomorphic to
$ X\mtim (\bA^{1}\backslash \{0\}) $.
Therefore $ p $ gives a deformation of
$ T_{Z}X $ to
$ X $,
see [36].

Let
$ M $ be a regular holonomic
$ \cD_{X} $-module.
Let
$$
\tM = \mopls_{i\in \bZ} M\otimes t^{i}.
$$
This has naturally a structure of
$ \cD_{X}\otimes_{\bC}\bC[t,t^{-1}]\langle\partial_{t}
\rangle $-module, and is identified with the pull-back of
$ M $ by the projection
$ q : \tX^{*} \to X $.
Viewed as a
$ \cD_{X}\otimes_{\bC}\bC[t]\langle\partial_{t}\rangle $-module,
$ \tM $ is identified with
$ \rho_{\ssbull}j_{*}q^{*}M $ where
$ j : \tX^{*} \to \tX $ and
$ \rho : \tX\to X $ are natural morphisms.
(Here
$ \rho_{\ssbull} $ and
$ j_{*} $ are direct images as Zariski sheaves.
Note that the direct image of
$ \cD $-modules
$ j_{*} $ is defined by the sheaf-theoretic direct image
in the case of open embeddings, see [5].)

Consider the filtrations
$ V $ of Kashiwara and Malgrange on
$ M $ along
$ Z $,
and on
$ j_{*}q^{*}M $ along
$ T_{Z}X = p^{-1}(0) $.
It is known that we have canonical isomorphisms
$$
\rho_{\ssbull}V^{\alpha-r+1}(j_{*}q^{*}M) =
V^{\alpha}\tM := \mopls_{i\in \bZ} V^{\alpha -i}
M\otimes t^{i}.
\leqno(1.3.1)
$$
In particular,
$$
\rho_{\ssbull}\Gr_{V}^{\alpha-r+1}(j_{*}q^{*}M) =
\Gr_{V}^{\alpha}\tM = \mopls_{i\in \bZ}\Gr_{V}^{\alpha -i}
M\otimes t^{i}.
$$
(These are crucial to relate Kashiwara's construction [16] with
the Verdier specialization [36].)
Here the shift of the filtrations
$ V $ comes from the difference of the codimensions.
The above assertion can be verified using the
$ \bC^{*} $-action on
$ \tX $ (or the action of the corresponding vector field)
which comes from the natural
$ \bC^{*} $-action on
$ \tX^{*} = X\mtim (\bA^{1}\backslash \{0\}) $ and which
corresponds to the grading by the order of
$ t $.
Indeed, if
$ (x_{1},\dots,x_{n}) $ is a local coordinate system of
$ X $ such that
$ Z $ is locally given by
$ x_{i}=0 $ for
$ i\le r $, then
$ \tX $ has a local coordinate system
$ (\tx_{1},\dots,\tx_{n},\tilt) $ such that
$ \tx_{i} = x_{i}/t \,(i \le r) $,
$ \tx_{i} = x_{i} \,(i > r) $,
$ \tilt = t $ on
$ \tX^{*} $, and hence the vector field corresponding to
the
$ \bC^{*} $-action is given by
$$
t\frac{\partial}{\partial t} =
\tilt\frac{\partial}{\partial\tilt} -
\sum_{i\le r}\tx_{i}\frac{\partial}{\partial\tx_{i}}.
\leqno(1.3.2)
$$
Note that
$ M\otimes t^{i} $ is identified with
$ \Ker(t\partial/\partial t-i) \subset \tM $ so that
$ \tilt\partial/\partial \tilt - i $ is identified with
$ \sum_{i\le r}\tx_{i}\partial/\partial \tx_{i} $ on
$ M\otimes t^{i} $.
Actually, the existence of the filtration
$ V $ can be reduced to the hypersurface case using this
argument.
The finite generatedness of (1.3.1) is
related to condition (iii) in (1.1) and to the property in Remark
(1.2) (ii).

We have also a canonical isomorphism
$$
\rho_{\ssbull}V^{k}\cD_{\tX} = \mopls_{i\in\bZ}
V^{k-i}\cD_{X}\otimes_{\bC}\bC[t\partial_{t}]t^{i}\quad
\hbox{for}\,\, k \ge 0,
\leqno(1.3.3)
$$
using
$ \partial/\partial\tx_{i} = t\partial/\partial x_{i} $,
$ \tx_{i} = t^{-1}x_{i} \,(i \le r) $, etc.
Then
$$
\rho_{\ssbull}\Gr_{V}^{k}\cD_{\tX} = \mopls_{i\in\bZ}
\Gr_{V}^{k-i}\cD_{X}\otimes_{\bC}\bC[t\partial_{t}]t^{i}\quad
\hbox{for}\,\, k \ge 0.
$$

The specialization
$ \SP_{Z}M $ of
$ M $ along
$ Z $ is defined by
$$
\SP_{Z}M = \psi_{t}(j_{*}q^{*}M) :=
\mopls_{0<\alpha\le 1}\Gr_{V}^{\alpha}(j_{*}q^{*}M),
$$
and its direct image by
$ \rho $ is identified with
$$
\mopls_{r-1<\alpha\le r}\Gr_{V}^{\alpha}\tM =
\mopls_{r-1<\alpha\le r}\mopls_{i\in \bZ}\Gr_{V}^{\alpha -i}
M\otimes t^{i}.
$$
(Here the shift of the indices comes from the difference of
the codimensions as above.)
For
$ \lambda = \exp(-2\pi i\alpha) $, we define the
$ \lambda $-eigen part by
$$
\SP_{Z,\lambda}M = \psi_{t,\lambda}(j_{*}q^{*}M) :=
\Gr_{V}^{\alpha}(j_{*}q^{*}M),
$$
which is identified with
$$
\Gr_{V}^{\alpha}\tM =
\mopls_{i\in \bZ}\Gr_{V}^{\alpha -i}M\otimes t^{i}.
$$

By the Riemann-Hilbert correspondence,
$ \SP_{Z}M $ corresponds via the analytic
de Rham functor
$ \DR $ to the specialization
$ \SP_{Z}K := \psi_{t}(\boR j_{*}q^{*}K) $ (see [10], [36]) of
$ K = \DR(M^{\rm an}) $ where
$ M^{\rm an} $ denotes the associated analytic
$ \cD $-module, see [16].
Note that
$ \SP_{Z,\lambda}M $ corresponds to
$ \SP_{Z,\lambda}K := \psi_{t,\lambda}(\boR j_{*}q^{*}K) $,
where
$ \psi_{t,\lambda} $ denotes the
$ \lambda $-eigen part of
$ \psi_{t} $.

\bigskip\bigskip
\centerline{{\bf 2. Bernstein-Sato polynomial}}

\bigskip\noindent
In this section we define the Bernstein-Sato polynomial for an
arbitrary variety (which is not necessarily reduced
nor irreducible), and prove its existence and
well-definedness together with Theorem~5.

\medskip\noindent
{\bf 2.1.
$ b $-Function.}
With the notation and the assumptions of (1.1), let
$ M $ be a quasi-unipotent regular holonomic
$ \cD_{X} $-module.
For a (local) section
$ m $ of
$ M $,
the {\it Bernstein-Sato polynomial} (i.e. the
$ b $-{\it function})
$ b_{m}(s) $ along
$ Z $ is defined to be the monic minimal polynomial of the
action of
$ s := -\theta -r $ on
$$
\oM_{m} := (V^{0}\cD_{X})m/(V^{1}\cD_{X})m.
\leqno(2.1.1)
$$

The action of
$ \theta $ on
$ \oM_{m} $ is independent of the choice of
$ \theta $.
The existence of
$ b_{m}(s) $ easily follows from that of the filtration
$ V $ of Kashiwara and Malgrange on
$ M $ along
$ Z $.
Indeed, the existence is equivalent to the finiteness of the induced
filtration
$ V $ on
$ \oM_{m} $, and setting
$ M' = \cD_{X}m $, it is sufficient to show for
$ \beta \gg 0 $
$$
V^{\beta}M' \subset (V^{1}\cD_{X})m,\quad
(V^{0}\cD_{X})m \subset V^{-\beta}M'.
\leqno(2.1.2)
$$
Since the induced filtration
$ V $ on
$ M' $ satisfies the conditions of the
$ V $-filtration (see Remark (iii) in (1.2)), we have for some
$ \beta_{0} $,
$ V^{\beta_{0}+i}M' = V^{i}\cD_{X}V^{\beta_{0}}M' $ for any
$ i > 0 $, and
$ V^{\beta_{0}}M' \subset (V^{j}\cD_{X})m $ for some
$ j $ by the coherence of
$ V^{\beta_{0}}M' $.
So the first inclusion follows.
The second inclusion is clear.

Note that
$ \alpha $ is a root of
$ b_{m}(-s) $ if and only if
$ \Gr_{V}^{\alpha}\oM_{m} \ne 0 $.
This implies
$$
\max\{\alpha : m\in V^{\alpha}M \} =
\min\{\alpha : b_{m}(-\alpha) = 0 \},
\leqno(2.1.3)
$$
because the left-hand side coincides with
$ \min\{\alpha : \Gr_{V}^{\alpha}((V^{0}\cD_{X})m)
\ne 0 \} $, which is strictly smaller than
$ \min\{\alpha : \Gr_{V}^{\alpha}((V^{1}\cD_{X})m)
\ne 0 \} $, see [26] for the codimension
$ 1 $ case.

\medskip\noindent
{\bf 2.2.~Proposition.}
{\it Let
$ i : X \to Y $ be a closed embedding of smooth varieties.
Let
$ Z $ be a smooth closed subvariety of
$ X $.
Let
$ M $ be a regular holonomic
$ \cD_{X} $-module, and
$ m $ be a section of
$ M $.
Let $i_{*}M$ and $i_{\ssbull}M$ denote the
direct images of $M$ as $\cD_{X}$-module and as $\cO_{X}$-module,
respectively.
Let
$ (y_{1}, \dots, y_{n}) $ be a local coordinate system of
$ Y $ such that
$ X = \{y_{i} = 0 \,(i\le q)\} $.
Let
$ m' $ be the element of
$ i_{*}M $ corresponding to
$ m\otimes 1 $ by the isomorphism
$ i_{*}M \simeq i_{\ssbull}M\otimes_{\bC}\bC[\partial_{1},
\dots, \partial_{q}] $ where
$ \partial_{i} = \partial /\partial y_{i} $, see
{\rm [5]}.
Then the Bernstein-Sato polynomial
$ b_{m}(s) $ of
$ m $ along
$ Z $ coincides with
$ b_{m'}(s) $ of
$ m' $ along
$ Z $, and
$$
V^{\alpha}(i_{*}M) = \msum_{\nu}i_{\ssbull}
V^{\alpha+|\nu|}M\otimes \partial^{\nu},
\leqno(2.2.1)
$$
where
$ V $ is the filtration of Kashiwara and Malgrange along
$ Z $,
$ \partial^{\nu} = \mprod_{i}\partial_{i}^{\nu_{i}} $ for
$ \nu = (\nu_{1}, \dots, \nu_{r}) \in \bN^{r} $, and
$ |\nu| = \msum_{i} \nu_{i} $.
}

\medskip\noindent
{\it Proof.}
Since the assertion is local, we may assume
$ \codim_{Y}X = 1 $,
and furthermore
$ Z = \{y_{i} = 0 \,(1 \le i \le r+1)\} $ in
$ Y $.
Then
$ \theta $ on
$ X $ and
$ Y $ can be given respectively by
$ \theta_{X} := \msum_{2\le i\le r+1}y_{i}\partial_{i} $
and
$ \theta_{Y} := \msum_{1\le i\le r+1}y_{i}\partial_{i} $.
Note that
$$
\theta_{X}+r = \msum_{2\le i\le r+1}\partial_{i}y_{i},\quad
\theta_{Y}+r+1 = \msum_{1\le i\le r+1}\partial_{i}y_{i}.
$$
Since
$ m' $ is annihilated by
$ y_{1} $,
we see that
$$
V^{k}\cD_{Y}\,m' = \mopls_{i\ge 0}(V^{i+k}\cD_{X}\,m)
\otimes \partial_{1}^{i}\quad \text{for}\,\, k = 0, 1.
\leqno(2.2.2)
$$
As
$ V^{i}\cD_{X}\,m/V^{i+1}\cD_{X}\,m $ is annihilated by
$ b_{m}(s+i) $,
the first assertion follows.
The proof of the second assertion is similar.

\medskip\noindent
{\bf 2.3.~Remark.}
With the above notation, assume
$ M $ has the Hodge filtration
$ F $.
Then the Hodge filtration
$ F $ on
$ i_{*}M $ is given by
$$
F_{p}(i_{*}M) = \msum_{\nu}
i_{\ssbull}F_{p-|\nu|}M\otimes \partial^{\nu}.
\leqno(2.3.1)
$$
In particular,
$ F_{p_{0}}(i_{*}M) = i_{\ssbull}F_{p_{0}}M $ locally if
$ p_{0} = \min\{p : \Gr_{p}^{F}M\ne 0 \} $.
Globally, we need the twist by the relative dualizing sheaf.

If
$ M = \cO_{X} $, we have
$ \Gr_{p}^{F}\cO_{X} = 0 $ for
$ p \ne -n $ where
$ n = \dim X $.
In particular.
$ p_{0} = -n $.

\medskip\noindent
{\bf 2.4.~The graph embedding.}
Let
$ X $ be a smooth affine variety, and let
$ Z $ be a (not necessarily reduced nor irreducible) closed
subvariety with
$ f_{1}, \dots, f_{r} $ generators of the ideal of
$ Z $ in
$ X $ (where
$ f_{j} \ne 0 $ for any
$ j $).
Let
$ i_{f} : X \to X' := X\mtim \bA^{r} $ be the graph embedding
of
$ f := (f_{1}, \dots, f_{r}) : X \to \bA^{r} $, i.e.
$ i_{f}(x) = (x,f_{1}(x), \dots, f_{r}(x)) $.
For a quasi-unipotent regular holonomic
$ \cD_{X} $-module
$ M $, let
$ M' = (i_{f})_{*}M $ be the direct image as
$ \cD $-module.
We have a natural isomorphism
$$
M' = (i_{f})_{\ssbull}M\otimes_{\bC}\bC[\partial_{1},
\dots, \partial_{r}],
\leqno(2.4.1)
$$
where
$ (i_{f})_{\ssbull} $ denotes the sheaf-theoretic direct
image, and
$ \partial_{i} = \partial /\partial t_{i} $ with
$ (t_{1}, \dots, t_{r}) $ the canonical coordinates of
$ \bA^{r} $.
Furthermore, the action of
$ \cO_{X}[\partial_{1}, \dots, \partial_{r}] $ on
$ M' $ is given by the canonical one (without using
$ f) $,
and the action of a vector field
$ \xi $ on
$ X $ and that of
$ t_{i} $ are given by
$$
\aligned
\xi (m\otimes \partial^{\nu})
&= \xi m\otimes \partial^{\nu} -
\msum_{i} (\xi f_{i})m\otimes \partial^{\nu +1_{i}},
\\
t_{i}(m\otimes \partial^{\nu})
&= f_{i}m\otimes \partial^{\nu} -
\nu_{i}m\otimes \partial^{\nu-1_{i}},
\endaligned
$$
where
$ \partial^{\nu} = \mprod_{i} \partial_{i}^{{\nu}_{i}} $
with
$ \nu = (\nu_{1}, \dots, \nu_{r}) \in \bN^{r} $,
and
$ 1_{i} $ is the element of
$ \bZ^{r} $ whose
$ j $-th component is
$ 1 $ if
$ j = i $ and
$ 0 $ otherwise.
In the case
$ M = \cO_{X} $, we have a canonical injection
$$
\cO_{X}\otimes_{\bC}\bC[\partial_{1}, \dots, \partial_{r}]
\hookrightarrow
\cO_{X}[\mprod_{i}f_{i}^{-1},s_{1},\dots, s_{r}]
\mprod_{i}f_{i}^{{s}_{i}},
\leqno(2.4.2)
$$
such that
$ s_{i} $ is identified with
$ -\partial_{i}t_{i} $, see [15], [21], [27].

Let
$ V $ be the filtration of Kashiwara and Malgrange on
$ \cD_{X'} $ along
$ X\mtim \{0\} $.
Let
$ \theta = \msum_{i} t_{i}\partial_{i} $ and
$ \theta^{*} = -\msum_{i} \partial_{i}t_{i} \,(= -\theta-r) $.
(Actually
$ {}^{*} $ comes from the involution of the ring of differential
operators, which is used in the transformation between left and
right
$ \cD $-modules.)
For a (local) section
$ m $ of
$ M $,
the {\it Bernstein-Sato polynomial} (i.e. the
$ b $-{\it function})
$ b_{f,m}(s) $ is defined to be that for
$ m\otimes 1 $.
This is the minimal polynomial of the action of
$ \theta^{*} $ on
$$
\oM_{f,m} := V^{0}\cD_{X'}(m\otimes 1)/V^{1}
\cD_{X'}(m\otimes 1).
\leqno(2.4.3)
$$
Note that the roots of
$ b_{f,m}(s) $ are rational numbers because the filtration
$ V $ is indexed by rational numbers.

If
$ M = \cO_{X} $ and
$ m = 1 $, then
$ b_{f,m}(s) $ is denoted by
$ b_{f}(s) $, and
$ b_{Z}(s) = b_{f}(s-r') $ with
$ r' = \codim_{X}Z $.
This definition of the Bernstein-Sato polynomial coincides
with the one in the introduction by (1.1.1) and (2.4.2), because
$ \theta^{*} $ belongs to the center of
$ \Gr_{V}^{0}\cD_{X'} $.
(Indeed, if
$ \sum_{j\le r}\mu_{j} - \sum_{i\le r}\nu_{i} = 1 $ with
$ \mu_{j},\nu_{i}\in\bN $, then
$ \mu_{j} \ge 1 $ for some
$ j $, and (1.1.1) for
$ k = 0 $ is contained in the
$ \bC $-algebra generated by
$ x_{i}\partial_{x_{j}} \,(i,j\le r) $ and
$ \partial_{x_{j}}\,(j>r) $.)

If
$ Z $ is not affine, then
$ b_{Z}(s) $ is defined to be the least common
multiple of
$ b_{(Z,z)}(s - \dim(Z,z) + \dim Z) $ for
$ z \in Z $, where
$ b_{(Z,z)}(s) $ is the
$ b $-function of a sufficiently small affine
neighborhood of
$ z $ in
$ Z $.

\medskip\noindent
{\bf 2.5.~Theorem.}
{\it The Bernstein-Sato polynomial
$ b_{f}(s) $ is independent of the choice of
$ f = (f_{1}, \dots, f_{r}) $ {\rm (}provided that
$ \dim X $ is fixed\,{\rm )}, and
$ b_{Z}(s) $ depends only on
$ Z $.
}

\medskip\noindent
{\it Proof.}
We first show that
$ b_{f}(s) $ is independent of the choice of the
generators
$ f_{1}, \dots, f_{r} $ if
$ X $ is fixed.
Let
$ g_{1}, \dots, g_{r'} $ be other generators.
Then we have
$ g_{i} = \msum_{j} a_{i,j}f_{j} $ with
$ a_{i,j} \in \cO_{X} $.
Set
$ f = (f_{1}, \dots, f_{r}) $,
$ g = (g_{1}, \dots, g_{r'}) $.
Let
$ i_{f} : X \to X\mtim \bA^{r} $,
$ i_{f,g} : X \to X\mtim \bA^{r+r'} $ be the embeddings by
the graphs of
$ f $ and
$ (f,g) $.
Consider an embedding
$ \phi : X\mtim \bA^{r} \to X\mtim \bA^{r+r'} $
defined by
$ \phi (x,s_{1}, \dots, s_{r}) = (x,s_{1}, \dots, s_{r},
s'_{1}, \dots, s'_{r'}) $ with
$ s'_{i} = \msum_{j} a_{i,j}s_{j} $.
Then
$ \phi\scirc i_{f} = i_{f,g} $.
So the independence of the choice of
$ f_{1}, \dots, f_{r} $ with
$ X $ fixed follows from Proposition (2.2).

Now we have to show
$ b_{g}(s) = b_{f}(s+1) $ for
$ g = (f_{1}, \dots, f_{r}, x) $ on
$ Y := X\mtim\bA^{1} $, where
$ x $ is the coordinate of
$ \bA^{1} $ and
$ f_{1}, \dots, f_{r}, x $ are viewed as generators of
the ideal of
$ Z\mtim\{0\} $ in
$ X\mtim\bA^{1} $.
Since the
$ b $-function
$ b_{x}(s) $ of
$ x $ is
$ s+1 $ as well-known,
this is a special case of Theorem~5, and
follows from (2.9) below.
Since the construction in (2.4) is compatible with
the pull-back by an etale morphism,
the assertion in the Proposition follows from Proposition (2.2).

\medskip\noindent
{\bf 2.6. Proposition.}
{\it With the notation of {\rm (2.4)}, let
$ V $ be the filtration on
$ M $ induced by the
$ V $-filtration on
$ M' $ along
$ X\mtim\{0\} $ using the isomorphism {\rm (2.4.1)}
and identifying
$ M $ with
$ (i_{f})_{\ssbull}M\otimes 1 $.
Then the filtration
$ V $ on
$ M $ is independent of the choice of
$ f_{1},\dots,f_{r} $.
}

\medskip\noindent
{\it Proof.}
This follows from (2.2.1) using the same argument as
in (2.5).

\medskip\noindent
{\bf 2.7.~Proposition.}
{\it With the notation of {\rm (1.3)}, let
$ m \in M $, and set
$ \tm = m\otimes 1 \in \tM $.
Let
$ b_{m}(s) $ and
$ b_{\tm}(s) $ be the
$ b $-functions of
$ m $ and
$ \tm $ along
$ Z $ and
$ T_{Z}X = p^{-1}(0) $ respectively.
Then
$ b_{\tm}(s) = b_{m}(s+r-1) $.
}

\medskip\noindent
{\it Proof.}
This follows from (1.3.2) and (1.3.3).

\medskip\noindent
{\bf 2.8.~Corollary.}
{\it With the notation of {\rm (2.7)}, let
$ \alpha_{j} $ be the roots of
$ b_{m}(-s) $.
Then the
$ \exp(2\pi i\alpha_{j}) $ are eigenvalues of the monodromy on
the nearby cycle sheaf
$ \psi_{t}(\boR j_{*}q^{*}K) $ in the notation of {\rm (1.3)}.
Furthermore, if
$ m $ generates
$ M $, then the
$ \exp(2\pi i\alpha_{j}) $ coincide with the eigenvalues of the
monodromy on the nearby cycle sheaf.
}

\medskip\noindent
{\it Proof.}
This follows from Proposition (2.7) together with
(1.3) and [16], [22].

\medskip\noindent
{\bf 2.9. Proof of Theorem~5.}
Let
$ Z = X\times Y $ and
$ Z' = Z\times \bA^{r+r'} $ (similarly for
$ X', Y' $ with
$ Z, r+r' $ replaced by
$ X,r $ and
$ Y,r' $ respectively).
Let
$ M' = (i_{f})_{*}\cO_{X},
N' = (i_{g})_{*}\cO_{Y},
R' = (i_{h})_{*}\cO_{Z} $.
It is clear that
$ R' $ is the external product
$ M' \boxtimes N' \,
(:= pr_{1}^{*}M' \otimes_{\cO} pr_{2}^{*}N') $ of
$ M' $ and
$ N' $.

Define a filtration
$ G $ on
$ M', N', R' $ by
$ G^{i}M' = V^{i}\cD_{X'}(1\otimes 1) $,
and similarly for $N'$ and $R'$.
Then
$$
\aligned
G^{k}R'
&= \msum_{i+j=k} G^{i}M' \boxtimes G^{j}N',
\\
\Gr_{G}^{0}R
&= \mopls_{i+j=0} \Gr_{G}^{i}M' \boxtimes \Gr_{G}^{j}N'.
\endaligned
$$
Put
$ b'_{h}(s) = \prod_{\gamma}(s+\gamma)^{q_{\gamma}} $ with
$ q_{\gamma} $ as in Theorem~5.
Since
$ \Gr_{G}^{i}M' $ is annihilated by
$ b_{f}(s+i) $ with
$ s = \msum_{1\le j\le r} s_{j} $,
and similarly for
$ \Gr_{G}^{j}N' $, we see that
$ \Gr_{G}^{i}M' \boxtimes \Gr_{G}^{j}N' $ is annihilated by
$ b'_{h}(s+i+j) $ with
$ s = \sum_{1\le j\le r+r'} s_{j} $.
Thus
$ b_{h}(s) $ divides
$ b'_{h}(s) $.
Moreover, we get the equality
$ b_{h}(s) = b'_{h}(s) $ by looking at the action of
$ \sum_{1\le j\le r+r'} s_{j} $ on
$ \Gr_{G}^{0}M' \boxtimes \Gr_{G}^{0}N' $.
This completes the proof of Theorem~5.

\medskip\noindent
{\bf 2.10.~Another description of the Bernstein-Sato polynomial.}
For
$ c = (c_{1}, \dots, c_{r}) \in \bZ^{r} $, let
$ I(c)_{-} = \{i : c_{i} < 0 \} $.
Then the Bernstein-Sato polynomial
$ b_{f}(s) $ is the monic polynomial of the smallest degree
such that
$ b_{f}(s)\mprod_{i}f^{s_{i}} $ belongs to the
$ \cD_{X}[s_{1}, \dots, s_{r}] $-submodule generated by
$$
\mprod_{i\in I(c)_{-}}\hbox{$\binom{s_{i}}{-c_{i}}$}\cdot
\mprod_{i=1}^{r}f_{i}^{s_{i}+c_{i}},
\leqno(2.10.1)
$$
where
$ c = (c_{1}, \dots, c_{r}) $ runs over the elements of
$ \bZ^{r} $ such that
$ \msum_{i}c_{i} = 1 $.
Here
$ s = \msum_{i=1}^{r}s_{i} $, and
$ \binom{s_{i}}{m} = s_{i}(s_{i}-1)\cdots (s_{i}-m+1)/m! $
as usual.

This definition of the Bernstein-Sato polynomial coincides with
the ones in the introduction and in (2.4).
Indeed,
$ s_{i}t_{i}^{-1} $ corresponds to
$ -\partial_{t_{i}} $ and the relation
$ t_{i}^{-1}s_{i} = (s_{i}-1)t_{i}^{-1} $ implies
$$
(s_{i}t_{i}^{-1})^{-c_{i}} = (-c_{i})!\hbox{$\binom{s_{i}}{-c_{i}}$}
t^{c_{i}} \quad\text{for}\quad c_{i} < 0.
$$
Then we put
$ \theta_{i} = t_{i} $ if
$ c_{i} > 0 $, and
$ \theta_{i} = \partial_{t_{i}}^{-1} $ if
$ c_{i} < 0 $, and consider
$ \prod_{i}\theta_{i}^{c_{i}} $ for
$ c \in \bZ^{r} $ with
$ \sum_{i}c_{i} = 1 $.

\newpage
\centerline{\bf 3. Proofs of Theorems 1--4}

\bigskip\noindent
In this section, we prove Theorems 1--4
using the bistrictness of the direct image.
We first recall the definition of multiplier ideals, see [20], [25].

\medskip\noindent
{\bf 3.1.~Multiplier ideals.}
Let
$ X $ be a smooth variety, and
$ Z $ be a closed subvariety of
$ X $ which is not necessarily reduced
nor irreducible.
Let
$ \pi : Y \to X $ be an embedded resolution of
$ Z $, i.e.
$ Y $ is smooth and
$ D :=\pi^{-1}(Z) $ is a divisor with normal crossings.
Here the ideal of
$ D $ is generated by the pull-back of the ideal of
$ Z $, and we assume it is locally principal.
Let
$ D_{i} $ be the irreducible components of
$ D $ with multiplicity
$ m_{i} $.
Let
$ \omega_{Y/X}\ $ be the relative dualizing sheaf
$ \omega_{Y}\otimes_{\pi^{-1}\cO_{X}}\pi^{-1}\omega_{X}^{\vee} $
(where
$ \omega_{X}^{\vee} $ is the dual line bundle of
$ \omega_{X} $).
Then, for a positive rational number
$ \alpha $, the multiplier ideal
$ \cJ(X,\alpha Z) $ is defined by
$$
\cJ(X,\alpha Z) = \pi_{\ssbull}(\omega_{Y/X}\otimes_{\cO_{Y}}
\cO_{Y}
(-\msum_{i}[\alpha m_{i}]D_{i})),
\leqno(3.1.1)
$$
using the trace morphism
$ \pi_{\ssbull}\omega_{Y/X}\to \cO_{X} $.
They define a decreasing filtration on
$ \cO_{X} $ which is indexed discretely and
right-continuously.
Note that
$ V $ is left-continuously indexed, see (1.1).

\medskip
The following is the key to the proof of Theorem~1.

\medskip\noindent
{\bf 3.2.~Proposition.}
{\it Let
$ X, Z $ be as in {\rm (1.1)}, and
$ Y $ be a smooth projective variety.
Let
$ \pr : X\mtim Y \to X $ denote the first projection, and
$ (M,F) $ be a filtered
$ \cD $-module underlying a mixed Hodge module on
$ X\mtim Y $.
Let
$ V $ be the filtration of Kashiwara and Malgrange on
$ M $ along
$ Z\mtim Y $.
Put
$ p_{0} = \min\{p : \Gr_{p}^{F}M\ne 0 \} $.
Then the direct image of
$ V $ gives the filtration
$ V $ of Kashiwara and Malgrange along
$ Z $, and the bifiltered direct image is bistrict.
In particular, we have a natural isomorphism
$$
F_{p_{0}}V^{\alpha}\cH^{i} \pr_{*}M =
R^{i}\pr_{\ssbull}(\omega_{Y}\otimes F_{p_{0}}V^{\alpha} M),
\leqno(3.2.1)
$$
where
$ \pr_{*} $ denotes the direct image of bifiltered
$ \cD $-modules which is defined by using the relative
de Rham complex as in {\rm [6]}, and
$ R^{j}\pr_{\ssbull} $ is the sheaf-theoretic higher direct
image.
}

\medskip\noindent
{\it Proof.}
This is reduced to the case of codimension
$ 1 $ because the construction in (1.3) is compatible with
the direct image by the projection as above.
Indeed, the direct image of
$ (M;F,V) $ by
$ \pr $ is defined by using the sheaf-theoretic direct image of
the relative de Rham complex
$ \DR_{X\times Y/X}(M;F,V) $ whose
$ i $-th component is given by
$$
\Omega_{Y}^{\dim Y+i}\otimes_{\cO_{Y}}(M;F[-i],V),
$$
where the pull-back of
$ \Omega_{Y}^{\dim Y+i} $ by the second projection is omitted to
simplify the notation, and
$ (F[-i])_{p} = F_{p+i} $.
Furthermore, the Hodge filtration
$ F_{p_{0}-1} $ on
$ j_{*}q^{*}M $ is obtained by taking the intersection of
the direct image
$ j_{\ssbull}q^{*}F_{p_{0}}M $ with
$ V^{0} $ (see [28], 3.2.3), where the base changes of
$ j,q $ are also denoted by the same symbols,
and the shift of the Hodge filtration by
$ 1 $ comes from the smooth pull-back
$ q^{*} $.
Since
$ \rho_{\ssbull}j_{\ssbull}q^{*}F_{p_{0}}M $ is identified with
$ \mopls_{k\in\bZ}F_{p_{0}}M\otimes t^{k} $ (on which the action of
$ t $ is bijective), the above argument together with (1.3.1) implies
$$
\rho_{\ssbull}F_{p_{0}-1}j_{*}q^{*}M =
\mopls_{k\in\bZ}F_{p_{0}}V^{r-1-k}M\otimes t^{k}.
\leqno(3.2.2)
$$
Note that
$ F_{p_{0}}V^{r-1-k}M = F_{p_{0}}M $ for
$ k\gg 0 $.

By [28], 3.3.17, we have the bistrictness of
$ \tpr_{*}(j_{*}q^{*}M;F,V) $.
By (1.3.1) and (3.2.2), the filtrations
$ F $ and
$ V $ are compatible with the grading by the
powers of
$ t $ (i.e. with the
$ \bC^{*} $-action).
So the bistrictness of
$ \tpr_{*}(j_{*}q^{*}M;F,V) $ implies that of
$ pr_{*}(M;F,V) $ by the definition of direct image using
the relative de Rham complex,
where we use the direct factor of (3.2.2) for
$ k \gg 0 $.

Note that (3.2.1) follows from the bistrictness if
$ V $ on the left-hand side is the induced filtration on the direct
image
$ \cH^{i}pr_{*}M $, see also (3.3.1) below.
So it remains to show that the induced filtration
$ V $ on
$ \cH^{i}pr_{*}M $ coincides with the filtration of Kashiwara and
Malgrange.
But this is proved in [28], 3.3.17 for
$ \tpr_{*}(j_{*}q^{*}M;F,V) $, and the assertion for
$ pr_{*}(M;F,V) $ follows by using (1.3.1) for
$ M $ and
$ \cH^{i}pr_{*}M $.
This completes the proof of Proposition (3.2).

\medskip\noindent
{\bf 3.3.~Remark.}
Let
$ (K;F,V) $ be a bifiltered complex representing the direct
image
$ \pr_{*}(M;F,V) $ defined in the derived category of bifiltered
$ \cD $-modules.
Let
$ G $ be a finite filtration on
$ K^{j} $ defined by
$ G^{0}K^{j} = \Ker \,d $,
$ G^{1}K^{j} = \Im \,d $.
Then (3.2.1) is essentially equivalent to
$$
F_{p_{0}}V^{\alpha}\Gr_{G}^{0}K^{j} =
\Gr_{G}^{0}F_{p_{0}}V^{\alpha}K^{j},
\leqno(3.3.1)
$$
and it is nontrivial because of the problem of three
filtrations [9].
We have (3.3.1) if
$ (K;F,V) $ is bistrict, see [28], 1.2.13.

\medskip\noindent
{\bf 3.4.~Proof of Theorem~1.}
Let
$ \pi : Y \to X $ be a projective morphism of smooth varieties
such that
$ \pi^{-1}(Z) $ is a divisor with normal crossings on
$ Y $ as in (3.1), and
$ \pi $ induces an isomorphism over
$ X \setminus Z $.
Let
$ g = (g_{1}, \dots, g_{r}) $ with
$ g_{i} = f_{i}\scirc\pi $ so that
$ g = \pi^{*}f $.
Then
$ i_{f}\scirc\pi = (\pi\mtim id)\scirc i_{g} $.

Let $D$ be the effective divisor defined by the pull-back
of the ideal of $Z$, and let
$Y'={\mathcal Spec}_{\cO_{Y}}(\mopls_{m\geq 0}\cO_{Y}(-mD))$
be the line bundle over $Y$ corresponding to the invertible
sheaf
$ \cO_{Y}(D) $.
Let
$ i_{D} : Y \to Y' $ be the closed embedding
induced by $\cO_Y(-D)\hookrightarrow\cO_Y$.
As a section of a line bundle, this corresponds to
$ 1 \in \cO_{Y}(D) $.
We have an embedding $i : Y'\hookrightarrow Y\times {\mathbb A}^n$
induced by
the surjective morphism
$$
(g_{1}, \dots, g_{r}) : \mopls_{j}\cO_{Y}\to \cO(-D),
$$
by passing to symmetric algebras over $\cO_Y$.
Then
$ i\scirc i_{D} = i_{g} $, and we get a commutative diagram
$$
\CD
Y @>{i_{D}}>> Y' \\
@| @VV{i}V \\
Y @>{i_{g}}>> Y\times\bA^{r} \\
@V{\pi}VV @VV{\pi\times id}V \\
X @>{i_{f}}>> X\times\bA^{r} \\
\endCD
$$
Taking a local equation
$ g' $ of
$ D $, we have a local trivialization of the line bundle
$ Y'\to Y $, and the embedding
$ i_{D} : Y \to Y' $ is locally identified with the graph embedding
$ i_{g'} : Y \to Y\times\bA^{1} $.
Since
$ g'{}^{-1}(0) $ is a divisor with normal crossings,
the assertion of Theorem~1 for
$ g'{}^{-1}(0) \subset Y $ follows from [29], Prop. 3.5,
see also [6], Prop. 2.3.
Here the filtration
$ V $ on
$ \cO_{Y} $ is induced by the
$ V $-filtration along the zero section of the line bundle
$ Y'\to Y $ using the inclusion
$ i_{D} : Y \to Y' $ (which is locally identified with
$ i_{g'} : Y \to Y\times\bA^{1} $).

We have the factorization of
$ \pi : Y\to X $ by the closed embedding
$ i' : Y\to X\mtim Y $ and the projection
$ pr : X\mtim Y \to X $.
Since the
$ V $-filtration is compatible with the direct image under
a closed embedding by (2.2) and
$ F_{p_0}M $ does not essentially change by such direct images
(see (2.3)), it is enough to consider the direct image of
$ M := (i_{h}\scirc i')_{*}\cO_{Y} $ by
$ pr\mtim id $, where
$ h $ is the pull-back of
$ f $ by
$ pr $.
Here
$ p_{0} = -n $ with
$ n = \dim X $, see (2.3).

By Proposition (3.2) (applied to
$ M $ and
$ pr\mtim id $), it is then enough to show that
$$
F_{-n}V^{\alpha}\cH^{0}(pr\mtim id)_{*}M =
F_{-n}V^{\alpha}(i_{f})_{*}\cO_{X}
\quad\text{for}\,\,\, \alpha > 0.
\leqno(3.4.1)
$$
By the decomposition theorem for mixed Hodge modules [29],
$ ((i_{f})_{*}\cO_{X},F) $ is a direct factor of
$ \cH^{0}(pr\mtim id)_{*}(M,F) $, and the complement
$ (N',F) $ is supported on
$ X\mtim\{0\} $.
We have to show
$ V^{\alpha}N' = 0 $ for
$ \alpha > 0 $.
But any element of
$ N' $ is annihilated by a sufficiently high power of
the ideal
$ (t_{1}, \dots, t_{r}) $, and hence
it is annihilated by a polynomial
$ b(-s) $ in
$ -s $ whose roots are non positive integers.
(Note that it is annihilated by
$ -s $ if it is annihilated by
$ (t_{1}, \dots, t_{r}) $.)
So the assertion follows.

\medskip\noindent
{\bf 3.5.~Proof of Theorem~2.}
With the notation of (2.4), let
$ M = \cO_{X} $ and
$ n = \dim X $ so that
$ F_{-n}M' = \cO_{X} $.
Since the filtration
$ V $ on
$ \cO_{X} $ is the induced filtration, we have the injectivity
of
$$
\Gr_{V}^{\alpha}\cO_{X}\to
\Gr_{V}^{\alpha}((V^{0}\cD_{X'})(1\otimes 1)),
$$
and
$$
\min\{\alpha : \Gr_{V}^{\alpha}\cO_{X}\ne 0\} =
\min\{\alpha : \Gr_{V}^{\alpha}((V^{0}\cD_{X'})(1\otimes 1))
\ne 0\},
$$
which will be denoted by
$ \alpha'_{Z} $.
(Here the inequality
$ \le $ follows from the above injectivity, and we have the equality
because
$ (V^{0}\cD_{X'})m \subset
V^{\alpha}(i_{f})_{*}\cO_{X} $ for
$ m \in V^{\alpha}\cO_{X}\otimes 1 $.)
Then
$$
(V^{1}\cD_{X'})(1\otimes 1) \subset
V^{\alpha'_{Z}+1}(i_{f}){*}\cO_{X},
$$
and the minimal root of
$ b_{Z}(-s) $ coincides with
$ \alpha'_{Z} $.
Therefore the assertion follows from Theorem~1.

\medskip\noindent
{\bf 3.6.~Proof of Theorem~3.}
Since
$ Z $ is a local complete intersection and the assertion is
local, we may assume
$ \codim_{X} Z = r $ with the notation of (1.3).
Then
$ T_{Z}X $ is a trivial vector bundle
$ Z\mtim\bA^{r} $,
as
$$
\mopls_{i\ge 0} I_{Z}^{i}/I_{Z}^{i+1} \simeq
\cO_{Z}\otimes_{\bC} \bC[u_{1}, \dots, u_{r}]
$$
with
$ u_{i} $ corresponding to
$ f_{i} $.

Since
$ Z $ is reduced, the restriction of the specialization
$ \SP_{Z}\bQ_{X} $ to
$ T_{Z}X\mtim_{Z}Z_{\reg} \,(= Z_{\reg}\mtim\bA^{r}) $ is
a constant sheaf on it, where
$ Z_{\reg} $ is the largest smooth open subvariety of
$ Z $.
This implies that the intersection complex (see [2]) of
$ T_{Z}X $ (which is the pull-back of the intersection
complex of
$ Z $ up to a shift of complex) is a subquotient of
$ \SP_{Z,1}\bQ_{X} $, where
$ \SP_{Z,1} $ denotes the unipotent monodromy part of
$ \SP_{Z} $, see (1.3).

To apply the theory of
$ \cD $-modules, we take the embedding
$ i_{f} : X \to X' = X\mtim\bA^{r} $ defined by the graph of
$ f $, and consider the specialization along
$ X\mtim\{0\} $.
Let
$ M = (i_{f})_{*}\cO_{X} $ (as a direct image of a
$ \cD $-module).
We apply (1.3) to these.
Let
$ Z' := T_{Z}X = Z\mtim\bA^{r} $, and
$ M_{Z'} $ denote the
$ \cD_{X'} $-module corresponding to the intersection
complex of
$ Z' \subset X' $, where
$ X' $ is identified with the normal cone of
$ X\mtim\{0\} $ in
$ X' = X\mtim\bA^{r} $.
Let
$ \rho : X\mtim\bA^{r} \to X $ denote the projection.
Then the above argument implies that
$ M_{Z'} $ is a subquotient of
$ \Gr_{V}^{1}(j_{*}q^{*}M) $, i.e.
$ \rho_{\ssbull}M_{Z'} $ is a subquotient of
$ \Gr_{V}^{r}\tM $.

By [28], [29], these
$ \cD $-modules have Hodge filtrations, denoted by
$ F $.
Here we use the normalization as in the case of right
$ \cD $-module, see [6].
Let
$ n = \dim X $.
Then
$ \rho_{\ssbull}F_{p}M_{Z'} = F_{p}\Gr_{V}^{r}\tM = 0 $ for
$ p < -n $, and the restrictions of
$ \rho_{\ssbull}F_{-n}M_{Z'} $ and
$ F_{-n}\Gr_{V}^{r}\tM $ to
$ Z_{\reg} $ are both isomorphic to
the structure sheaf
$ \cO $ of
$ Z_{\reg} $ tensored by
$ \bC[u_{1},\dots,u_{r}] $ over
$ \bC $.
Note that
$ F_{-n}M = \cO_{X} $, and
$ F_{-n}M_{Z'} $ is isomorphic to the sheaf-theoretic direct image
of the relative dualizing sheaf of a resolution of singularities
of
$ Z' $, see [17], [30].
So we get
$$
\rho_{\ssbull}F_{-n}M_{Z'} = \omega_{X}^{\vee}\otimes_{\cO}\tom_{Z}
\otimes_{\bC}\bC[u_{1}, \dots, u_{r}],
\leqno(3.6.1)
$$
where
$ \omega_{X}^{\vee} $ is the dual of
$ \omega_{X} $.
Furthermore, the above argument implies that
$ \rho_{\ssbull}F_{-n}M_{Z'} $ is a subquotient of
$$
F_{-n}\Gr_{V}^{r}\tM =
\mopls_{i\in\bZ}\Gr_{V}^{r+i}\cO_{X}
\otimes t^{-i}
\leqno(3.6.2)
$$
in a compatible way with the grading induced by
the natural
$ \bC^{*} $-action.

We have
$ f_{i} \in V^{>r}\cO_{X} $, and
$ \cO_{X}/V^{>r}\cO_{X} = \cO_{Z} $ because
$ \cO_{Z} $ does not contain a submodule
whose support has strictly smaller dimension.
Since
$ \omega_{X}^{\vee}\otimes_{\cO}\tom_{Z} $ is a subquotient
of
$ \cO_{X}/V^{>r}\cO_{X} = \cO_{Z} $, it is a submodule by a
similar argument, and it is contained in
$ \Gr_{V}^{r}\cO_{X} $.
Then, defining
$ \cJ(X,rD)' $ to be the coherent ideal of
$ \cO_{X} $ corresponding to
$ \omega_{X}^{\vee}\otimes_{\cO}\tom_{Z} \subset
\Gr_{V}^{r}\cO_{X} $,
the assertion follows from Theorem~1.

\medskip\noindent
{\bf 3.7.~Proof of Theorem~4.}
Assume first that
$ Z $ has at most rational singularities.
By Theorems 2 and 3, we have
$ \alpha'_{f} = r $ because
$ -r $ is a root of
$ b_{f}(s) $.
Furthermore,
$ \tom_{Z} = \omega_{Z} $, and hence
$ \rho_{\ssbull}F_{-n}M_{Z'} = F_{-n}\Gr_{V}^{r}\tM $ in the above
notation.
(This is closely related with Corollary~1.)

Let
$ \tilt $ be as in (1.3.2), and define
$ N = \tilt\partial/\partial\tilt $ on
$ \Gr_{V}^{1}j_{*}q^{*}M $.
Since
$ M $ is pure of weight
$ n $,
$ q^{*}M $ is pure of weight
$ n + 1 $ by the property of external products,
see [29], 2.17.4.
(To simplify the notation, we do not shift the filtration
$ F $ as in [28], [29]
when we take
$ q^{*} $.
Note that the filtration is shifted to the opposite direction
when we take
$ \Gr_{1}^{V} $ in loc.~cit., and these two shifts cancel out.
Also we do not shift the complex as in [29] when we take
$ q^{*} $.
So
$ q^{*}M $ is a
$ \cD $-module.)

By definition ([28], 5.1.6) the weight filtration
$ W $ on
$ \Gr_{V}^{1}j_{*}q^{*}M $ is the monodromy filtration
shifted by
$ n $, i.e. it is characterized by the conditions:
$ NW_{i} \subset W_{i-2} $ and
$$
N^{j} : \Gr_{n+j}^{W}\simto\Gr_{n-j}^{W}\quad\text{for}
\quad j > 0.
\leqno(3.7.1)
$$
Furthermore,
$ \Gr_{n}^{W}\Gr_{V}^{1}j_{*}q^{*}M $ underlies a semisimple
Hodge module (see [28], 5.2.13), and
$ M_{Z'} $ is a direct factor of it.

Let
$$
K := \Ker(N :
\Gr_{V}^{1}j_{*}q^{*}M \to \Gr_{V}^{1}j_{*}q^{*}M).
$$
Then we have
$ M_{Z'} \subset \Gr_{n}^{W}K $ on
$ T_{Z}X\mtim Z_{\reg} $, and it holds everywhere by
the property of an intersection complex, see [2].
This implies that
$ F_{-n}\Gr_{V}^{r}M\otimes 1 $ is also contained in
$ \Gr_{n}^{W}K $, and hence the action of
$ \sum_{i\le r}\tx_{i}\partial/\partial \tx_{i} $ on
$ F_{-n}\Gr_{V}^{r}M\otimes 1 $ vanishes, see (1.3.2).
So the multiplicity of
$ -r $ as a root of the
$ b $-function is
$ 1 $, and
$ \alpha_{f} > r $.

Assume conversely
$ \alpha_{f} > r $, i.e.
$ \alpha'_{f} = r $ with multiplicity
$ 1 $ as a root of
$ b_{f}(-s) $.
Then
$ F_{-n}\Gr_{V}^{\alpha}M = 0 $ for
$ \alpha < r $, and
$ F_{-n}\Gr_{V}^{r}M = \cO_{Z} $, because
$ F_{-n}\Gr_{V}^{r}M $ is a quotient of
$ \cO_{Z} $ and there is no nontrivial
$ \cO_{Z} $-submodule of
$ \cO_{Z} $ supported on
$ \Sing\,Z $.
So it is enough to show that
$ \omega_{X}^{\vee}\otimes_{\cO}\tom_{Z} $ is a direct factor of
$ \cO_{Z} $, or more generally that
(3.6.1) is a direct factor of (3.6.2).

Since the multiplicity of
$ r $ is
$ 1 $, we see that
$ F_{-n}\Gr_{V}^{r}M \otimes 1 $ is contained in
$ K $, and so is (3.6.2) because
$ N $ is compatible with the action of
$ \mopls_{i\ge 0}I_{Z}^{i}/I_{Z}^{i+1}\otimes t^{-i} $.
(Here (3.6.2) can be calculated as in (3.8) below.)
Then, by the semisimplicity of the Hodge
module underlying
$ \Gr_{n}^{W}K $ (see [28], 5.2.13), it is sufficient to show
that
$ F_{-n}\Gr_{i}^{W}K = 0 $ for
$ i < n $.
(Here
$ \Gr_{i}^{W}K = 0 $ for
$ i > n $ by (3.7.1).)
But this follows from
$ F_{-n-1}\Gr_{V}^{1}j_{*}q^{*}M = 0 $ together with the fact
that
$ N $ is a morphism of type
$ (-1,-1) $, because the latter implies the isomorphisms
(see [28], 5.1.14)
$$
N^{j} : F_{p}\Gr_{n+j}^{W}\simto F_{p+j}\Gr_{n-j}^{W}\quad
\text{for}\quad j > 0.
$$
This completes the proof of Theorem~4.

\medskip\noindent
{\bf 3.8.~Proof of Corollary~1.}
We have
$ 1 \in V^{r}\cO_{X} $ by Theorem~3, and hence
$ \prod_{i}f_{i}^{\nu_{i}} \in V^{r+j}\cO_{X} $
with
$ j = |\nu| $.
Using the same argument as in (3.6),
we can show by induction on
$ j \ge 0 $ that
$ V^{>r+j}\cO_{X} $ is generated by
$ \prod_{i}f_{i}^{\nu_{i}} $ with
$ |\nu| = j + 1 $, and
$ \Gr_{V}^{r+j}\cO_{X} $ is a free
$ \cO_{Z} $-module generated by
$ \prod_{i}f_{i}^{\nu_{i}} $ with
$ |\nu| = j $ so that
$ \Gr_{V}^{\alpha}\cO_{X} = 0 $ for
$ \alpha \notin \bZ $.
Then the assertion follows from Theorem~1.

\medskip\noindent
{\bf 3.9.~Proof of Corollary~2.}
This follows from Theorem~1 using (2.1.3).

\medskip\noindent
{\bf 3.10.~Proposition.}
{\it With the notation and assumptions of {\rm (3.1)},
the union of the subgroups of
$ \bQ/\bZ $ generated by
$ 1/m_{j} $ contains the images of the roots of
$ b_{Z}(s) $ in
$ \bQ/\bZ $.
}

\medskip\noindent
{\it Proof.}
We may assume that the ideal of
$ Z $ is generated by
$ f_{1},\dots,f_{r} $.
By Corollary~(2.8), it is sufficient to consider the
eigenvalues of the monodromy on
$ \psi_{t}(\boR j_{*}q^{*}K) $, where
$ K = (i_{f})_{*}\bC_{X^{\rm an}} $.
Since the construction of the deformation space to the tangent cone
as in (3.1) is compatible with the embedded resolution
$ (Y,D) \to (X,Z) $, the assertion is reduced to the case where
$ Z $ is a divisor with normal crossings on a smooth variety,
using the commutativity of
the nearby cycle functor with the direct image under a
proper morphism together with the decomposition theorem [2].
In the divisor case, it is known that the eigen\-values of the
monodromy for
$ \psi_{t}(\boR j_{*}q^{*}K) $ coincide with those for
$ \psi_{g'}K $ where
$ g' $ is a local equation of the divisor.
(This follows also from (1.3.1).)
Furthermore, the eigen\-values of the monodromy in the normal
crossing case are
$ \exp(2\pi i k/m_{j}) \,(k \in \bZ) $ where the
$ m_{j} $ are the multiplicities of the irreducible components
of the divisor (see also [15], [22]).
So the assertion follows.

\bigskip\bigskip
\centerline{{\bf 4. Calculations in the monomial ideal case}}

\medskip\noindent
In this section we treat the case where the ideal of the subvariety
$ Z $ (which is not necessarily reduced nor irreducible) is a
monomial ideal, and calculate some examples.

\medskip\noindent
{\bf 4.1.~Monomial ideal case.}
Assume that
$ X $ is the affine space
$ \bA^{n} $ and the
$ f_{i} $ are monomials with respect to the coordinate
system
$ (x_{1},\dots,x_{n}) $ of
$ \bA^{n} $.
Write
$ f_{j} = \prod_{i=1}^{n}x_{i}^{a_{i,j}} $.
Let
$ \ell_{i}(\bs) = \sum_{j=1}^{r}a_{i,j}s_{j} $ for
$ \bs = (s_{1},\dots,s_{r}) $ so that
$$
\mprod_{j=1}^{r}f_{j}^{s_{j}} =
\mprod_{i=1}^{n}x_{i}^{\ell_{i}(\bs)}.
$$
Let
$ \ell(c) = (\ell_{1}(c),\dots,\ell_{n}(c)) $, and
$ I'(\ell(c))_{+} = \{i : \ell_{i}(c) > 0\} $.
Let
$ I(c)_{-} $ be as in (2.10), and define
$$
g_{c}(s_{1},\dots,s_{r}) =
\mprod_{j\in I(c)_{-}}\hbox{$\binom{s_{j}}{-c_{j}}$}\cdot
\mprod_{i\in I'(\ell(c))_{+}}
\hbox{$\binom{\ell_{i}(\bs)+\ell_{i}(c)}{\ell_{i}(c)}$}.
$$
Let
$ \fa_{f} $ be the ideal of
$ R := \bQ[s_{1},\dots,s_{r}] $ generated by
$ g_{c}(s_{1},\dots,s_{r}) $ where
$ c = (c_{1},\dots,c_{r}) $ runs over the elements of
$ \bZ^{r} $ such that
$ \msum_{i}c_{i} = 1 $ (see [32] for a more general case).
In a forthcoming paper [7], we will give a combinatorial description
of the roots of the Bernstein-Sato polynomial for monomial ideals.
For the proof, the following is used in an essential way.

\medskip\noindent
{\bf 4.2.~Proposition.}
{\it With the above notation and assumption,
the Bernstein-Sato polynomial
$ b_{f}(s) $ is the monic polynomial of smallest degree
such that
$ b_{f}(\sum_{i}s_{i}) $ belongs to the ideal
$ \fa_{f} $.
}

\medskip\noindent
{\it Proof.}
This follows from (2.10) using the
$ \bZ^{n} $-grading on
$ \cD_{X} $ such that the degree of
$ x_{i} $ is the
$ i $th unit vector
$ e_{i} $ of
$ \bZ^{n} $, and the degree of
$ \partial_{x_{i}} $ is
$ -e_{i} $.

\medskip\noindent
{\bf 4.3.~Finite generators of
$ \fa_{f} $.}
It follows from (4.2) that in order to compute the Bernstein-Sato
polynomial in the monomial case we need to solve the elimination problem
consisting in computing $\bQ[s_1+\ldots+s_r]\cap\fa_f$. This can be done
using a computer algebra program (e.g. \emph{Macaulay}2) if we can write
down
finitely many generators for the ideal $\fa_f$.
(Here we may also consider locally the subscheme of
$\bA^n$ defined by $\fa_f$,
and then take the direct image of its structure sheaf by the map
$ (s_{1},\dots,s_{r})\mapsto \sum_{i}s_{i} $, because
it is enough to determine the {\it annihilator} of the direct
image sheaf.
Using this, we can calculate the examples below by hand.)

We explain now how to find finite generators of
$ \fa_{f} $.
With the notation and the assumption of (4.1),
let
$ A_{\delta} = \{\sum_{i}c_{i}=\delta\}\subset \bZ^{r} $
for
$ \delta = 0, 1 $.
Let
$ \ell_{n+i}(c) = -c_{i} $ for
$ i = 1,\dots,r $.
For
$ \vep = (\vep_{1},\dots,\vep_{n+r})
\in \{1,-1\}^{n+r} $, define
$$
A_{\delta}^{\vep} = \{c\in
A_{\delta} : \vep_{i}\ell_{i}(c) \ge 0
\,\,\,\text{for}\,\,\,i = 1,\dots, n+r\}.
$$
If there is a subset
$ I(\vep) $ of
$ A_{1}^{\vep} $ for any
$ \vep $ such that
$$
A_{1}^{\vep} = \mcup_{c\in I(\vep)}(c+A_{0}^{\vep}),
\leqno(4.3.1)
$$
then we see that
$ \fa_{f} $ is generated by the
$ g_{c} $ for
$ c \in I(\vep) $ and
$ \vep \in \{1,-1\}^{n+r} $.
For each
$ \vep $ it is easy to show the existence of a finite subset
$ I(\vep) $ satisfying (4.3.1).
Probably some algorithms to find such finite subsets would
be already known.
(For example, if we find some
$ v\in A_{1}^{\vep} $, we can proceed by induction on
$ r $, considering the complement of
$ v + A_{0}^{\vep} $ and taking the restrictions of
$ A_{1}^{\vep} $ to the hyperplanes
$ \{\vep_{i}\ell_{i}(c) = m_{i}\} $ for
$ 0\le m_{i} < \vep_{i}\ell_{i}(v) $ and for every
$ i $.)

\medskip\noindent
{\bf 4.4.~Multiplier ideals of monomial ideals.}
With the above notation and assumption, let
$ \ba_{j} = (a_{1,j},\dots,a_{n,j}) $, and
$ P $ be the convex hull of
$ \mcup_{j}(\ba_{j}+\bR_{\ge 0}^{n}) $ in
$ \bR^{n} $.
Set
$ \bone = (1,\dots,1) \in \bZ^{n} $.
By [20] (see also [14]),
$ \cJ(X,\alpha Z) $ is generated by the monomials
$ x^{\nu} := \prod_{i}x_{i}^{\nu_{i}} $ with
$ \nu+\bone\in (\alpha+\vep)P $ where
$ \vep > 0 $ is sufficiently small.

As a corollary, the jumping coefficients are the numbers
$ \alpha $ such that
$ (\alpha\cdot\partial P) \cap \bZ_{>0}^{n} \ne \emptyset $
where
$ \partial P $ is the boundary of
$ P $.
Let
$ \phi_{j} $ be linear functions on
$ \bR^{n} $ such that the maximal dimensional faces of
$ P $ are contained in
$ \phi_{j} = 1 $.
Then the jumping coefficient corresponding to
$ x\in \bZ_{>0}^{n} $ (i.e. the number
$ \alpha $ such that
$ x\in \alpha\cdot\partial P $) is given by
$ \min\{\phi_{j}(x)\} $.

\medskip\noindent
{\bf 4.5.~Examples.} (i)
With the notation of (4.1), let
$ f_{i} = \mprod_{j\ne i}x_{j} $ for
$ i = 1,\dots, n $.
If
$ n = 3 $, then
$ b_{f}(s) = (s+3/2)(s+2)^{2} $.
Indeed, we have
$ \ell_{i}(\bs) = \sum_{j\ne i}s_{j} $, and
$ \fa_{f} $ is generated by
$ g_{(1,0,0)} = (\ell_{2}(\bs)+1)(\ell_{3}(\bs)+1) $ and
$ g_{(1,1,-1)} = s_{3}(\ell_{3}(\bs)+1)(\ell_{3}(\bs)+2) $
up to a permutation, see (4.3).
So the assertion follows by using
$ u_{i} := \ell_{i}(\bs)+1 $.
Note that if we consider only
$ g_{c} $ with all
$ c_{i} \ge 0 $, then we get only
$ (\ell_{2}(\bs)+1)(\ell_{3}(\bs)+1) $ up to a permutation,
and there is no nonzero polynomial
$ b(s) $ such that
$ b(\sum s_{i}) $ belongs to the ideal generated by these.

In this case, the polyhedron
$ P $ in (4.4) is defined by the inequalities
$ \sum_{i}x_{i}\ge 2 $ and
$ \sum_{i\ne j}x_{i}\geq 1 $ for all
$ j $.
Thus the jumping coefficients of this ideal are
$ k/2 $ for
$ k \ge 3 $ by (4.4), and
this is compatible with Theorem~2.

\medskip
(ii)
The above calculation can be extended to the case of a
monomial ideal generated by
$ x_{i}x_{j} $ for
$ 1 \le i < j \le n $ with
$ n > 3 $ (i.e.
$ Z $ is the union of coordinate axes).
In this case we have
$ b_{f}(s) = (s+\frac{n}{2})(s+\frac{n+1}{2})(s+n-1) $.
Note that
$ \frac{n}{2}, \frac{n+1}{2} $ are roots of
$ b_{f}(-s) $ by (4.4) together with Theorem~2, and
$ n-1 $ is a root of
$ b_{f}(-s) $ by restricting to a smooth point of
$ Z $.
(So it is enough to show that the above polynomial
belongs to
$ \fa_{f} $ to determine the
$ b $-function in this case, although it is not
difficult to generalize the above calculation.)

\medskip
(iii)
With the notation of (i), assume
$ n = 4 $.
In this case we can show that
$ b_{f}(s) = (s+4/3)(s+5/3)(s+3/2)(s+2)^{3} $.
Note that
$ 4/3 $ and
$ 3/2 $ are jumping coefficients, but
$ 5/3 $ is not.
(Indeed, the linear functions in (4.4) are given by
$ \sum_{k}x_{k}/3 $,
$ \sum_{k\ne i}x_{k}/2 $ and
$ \sum_{k\ne i,j}x_{k} $ for
$ i\ne j $.)

\medskip
(iv)
Assume
$ n = 3 $ and
$ f_{i} = x_{i}\prod_{j=1}^{3}x_{j} $ for
$ i = 1, 2, 3 $.
Then we have
$ b_{f}(s) = (s+3/4)(s+5/4)(s+6/4)(s+1)^{3} $.
Here
$ 3/4 $ and
$ 6/4 $ are jumping coefficients, but
$ 5/4 $ is not.
Note also that the support of
$ R/\fa_{f} $ is not discrete in this case.

\medskip\noindent
{\bf 4.6.~Integral closure of the ideal and
$ b $-function.}
It is known that the jumping coefficients depend
only on the integral closure of the ideal,
see [20], 11.1. For monomial ideals this follows also from
the description in (4.4).
However, this does not hold for the
$ b $-function.
For example, consider the ideal generated by
$ x^{2},y^{2} \in \bC[x,y] $.
Its
$ b $-function is
$ (s+1)(s+3/2)(s+2) $ by Theorem~4.
But the integral closure of this ideal is
generated by
$ x^{2},xy,y^{2} \in \bC[x,y] $, and one can check that its
$ b $-function is
$ (s+1)(s+3/2) $.

\newpage
\centerline{{\bf 5. Relation with other Bernstein-Sato polynomials}}

\bigskip\noindent
{\bf 5.1.~Bernstein-Sato polynomials of several variables.}
With the notation of (2.4),
let $w=(w_i)_i\in\bN^r$, and consider the ideal
$B_{w}(f)\subseteq\bC[s_1,\ldots,s_r]$
consisting of those polynomials $b$ satisfying
$$
b(s_1,\ldots,s_r)\prod_if_i^{s_i}\in\cD_X[s_i,\dots,s_r]\cdot
\prod_if_i^{s_i+w_i}.
$$
It was shown by Sabbah [27] that $B_{w}(f)\neq (0)$.

Suppose now that the $f_i$ are monomials as in (4.1). Then with the
notation of (4.1), we can show that the ideal $B_{w}(f)$ is generated by
$$b_{f,w}(s_1,\dots,s_r):=
\prod_{i=1}^n(\ell_i(\bs)+1)\cdots(\ell_i(\bs)+\ell_i(w)),$$
using the $ \bZ^{n} $-grading on $ \cD_{X} $ as in the proof of (4.2).
So this can be called the Bernstein-Sato polynomial of several variables
in this monomial case.
However, it does not seem easy to relate this with our Bernstein-Sato
polynomial of one variable. For example, if
$f_1=x^\alpha y,f_2=xy^\beta$ with $\alpha,\beta\ge 2$ and $n=r=2$,
then $$\ell_1=\alpha s_1+s_2,\quad\ell_2=s_1+\beta s_2,$$ and
the roots of the Bernstein-Sato polynomial of one variable $b_f(s)$
are $$-\frac{(\beta-1)i+(\alpha-1)j}{\alpha\beta-1}\quad
\text{for}\quad 1\le i\le \alpha,\,\, 1\le j\le \beta,$$
and $-1$ is also a root if it is not included there.

\medskip\noindent
{\bf 5.2.~Bernstein-Sato polynomials of one variable in [27, I].}
With the notation of (2.4),
$ b_{f}(s) $ coincides (up to a shift of variable) with the polynomial
$ b_{L} $ of the minimal degree satisfying the relation in
[27], I, 3.1.1 in the algebraic setting, if
$ \cM $ in loc.~cit. is the direct image of
$ \cO $ by the graph embedding, the multi-filtration
$ U $ is induced by that on
$ \cD_{X} $ using the action of
$ \cD_{X} $ on
$ 1\otimes 1 \in \cM $, and
$ L(s) = s_{1}+\cdots + s_{r} $, see also [13], 2.13.

If
$ X $ is affine space
$ \bA^{n} $, then
$ b_{f}(s) $ coincides with the
$ b $-function in [32], p.~194, if the weight vector is chosen
appropriately.
Algorithms to compute this
$ b $-function are given in loc.~cit., p.~196.

\medskip\noindent
{\bf 5.3.~Remark.}
The existence and the rationality of the roots of the polynomial
$ b_{\alpha} $ in [13], 2.13 can be reduced to the case
$ \alpha = (1,\dots,1) $ by considering the pull-back by the base
change of the finite covering
$ \pi : (\ttau_{j})\in \bC^{l} \mapsto (\tau_{j}) =
(\ttau_{j}^{\alpha_{j}}) \in \bC^{l} $ and taking the invariant by
the covering transformation group
$ \prod_{j}(\bZ/\alpha_{j}\bZ) $.
(Here we may assume that
$ \alpha \in \bZ_{>0}^{r} $ by replacing
$ r $ and changing the decomposition
$ X \times E $ in loc.~cit. if necessary.)
Then the assertions easily follow from the theory of Kashiwara and
Malgrange using (1.3.1).

\medskip\noindent
{\bf 5.4.~Bernstein-Sato polynomials of one variable in [27, II].}
With the notation of (2.4), we can define a polynomial
$ b'_{L} $ associated to a linear form
$ L : \bZ^{r}\to\bZ $ to be the monic polynomial of the smallest degree
satisfying the relation
$$
b'_{L}(L(\bs))\mprod_{i}f_{i}^{s_{i}} \in
\msum_{\nu} \cD_{X}[s_{1},\dots,s_{r}]\mprod_{i}f_{i}^{s_{i}+\nu_{i}},
\leqno(5.4.1)
$$
where the summation is taken over
$ \nu = (\nu_{1},\dots,\nu_{r}) \in \bZ^{r} $ such that
$ L(\nu) = 1 $, see [27], II, Prop. 1.1.
By (2.10), this polynomial
$ b'_{L} $ for
$ L(\bs) = \sum_{i}s_{i} $ divides our Bernstein-Sato polynomial
$ b_{f}(s) $, but they do not coincide,
and furthermore, Theorem~2 does not hold for
$ b'_{L} $ in general.
For example, consider the case
$ f_{1} = x^{3}, f_{2} = x^{2}y $ with
$ n = r = 2 $.
Then
$$
f_{1}^{s_{1}+1-k}f_{2}^{s_{2}+k} =
x^{\ell_{1}+3-k}y^{\ell_{2}+k} \quad\text{for}\quad 0\le k\le 3,
$$
where
$ \ell_{1} = 3 s_{1} + 2 s_{2}, \ell_{2} = s_{2} $.
Applying
$ \partial_{x}^{3-k}\partial_{y}^{k} $ to this, we see that
(5.4.1) holds with
$ b'_{L}(L(\bs)) $ replaced by
$$
\mprod_{1\le i\le 3-k}(\ell_{1}+i)\cdot\mprod_{1\le j\le k}(\ell_{2}+j).
$$
This implies
$ b'_{L}(s) = (s+2/3)(s+1)(s+4/3) $, and
$ -1/2 $ is not a root of it.
But this is not compatible with the jumping coefficients which include
$ 1/2 $ (see (4.4)),
and Theorem~2 does not hold for this polynomial.

Note also that the polynomial
$ b'_{L}(s) $ depends on the choice of generators of the ideal of
$ Z $.
Indeed,
if we put
$ g_{i} = f_{i} $ for
$ i \le r $ and
$ g_{r+1} = f_{1}^{2} $, then
$ g_{1}^{2}g_{r+1}^{-1} = 1 $ and the polynomial associated to
$ g_{1},\dots,g_{r+1} $ is
$ 1 $.

\end{document}